\documentclass{article}

\usepackage{latexsym,amssymb,graphicx,cite}
\input amssym.def \input amssym

\newtheorem{te}{Theorem}[section]

\newtheorem{de}[te]{Definition}
\newtheorem{lm}[te]{Lemma}

\newtheorem{pp}[te]{Proposition}
\newtheorem{co}[te]{Corollary}
\newtheorem{ex}[te]{Example}

\def\dokaz{\noindent{\bf Proof. }}
\def\kraj{\hfill $\Box$ \par \vspace*{2mm} }
\def\widemid{\hspace{1mm}\widetilde{\mid}\hspace{1mm}}

\newcommand{\zve}[1]{{{}^*\hspace{-0.5mm}#1}}
\newcommand{\zvez}[1]{{{}^*\hspace{-1mm}#1}}
\def\zvepar{\hspace{1mm}\zvez\parallel\hspace{1mm}}

\def\zvemid{\hspace{1mm}\zvez\mid\hspace{1mm}}
\def\nzvemid{\hspace{1mm}\zvez\nmid\hspace{1mm}}

\begin{document}

\begin{center}
           {\huge \bf A survey on divisibility of ultrafilters}\\
{\small \bf Boris  \v Sobot}\\[2mm]
{\small  Department of Mathematics and Informatics, University of Novi Sad,\\
Trg Dositeja Obradovi\'ca 4, 21000 Novi Sad, Serbia\\
e-mail: sobot@dmi.uns.ac.rs\\
ORCID: 0000-0002-4848-0678}
\end{center}

\begin{abstract}
An extension of the divisibility relation on $\mathbb{N}$ to the set $\beta\mathbb{N}$ of ultrafilters on $\mathbb{N}$ was defined and investigated in several papers during the last ten years. Here we make a survey of results obtained so far, adding several results connecting the themes of different stages of the research. The highlights include: separation of $\beta\mathbb{N}$ into the lower part $L$ (with its division into levels) and the upper part; identifying basic ingredients (powers of primes) and fragmentation of each ultrafilter into them; finding the corresponding upward closed sets belonging to an ultrafilter with given basic ingredients; existence and number of successors and predecessors of a given divisibility class; extending the congruence relation (in two ways) and checking properties of the obtained relations.
\end{abstract}



\section{Introduction}

Ultrafilters are a powerful tool, often used in mathematical logic, topology and other areas. Each ultrafilter establishes a kind of "majority criterion" which can be used in various constructions. Among important achievements of their usage are applications in infinitary combinatorics, which began with an ultrafilter proof of Hindman's theorem.

\begin{te}
For every coloring $f:\mathbb{N}\rightarrow\{1,2,\dots,k\}$ there is an infinite $A\subseteq\mathbb{N}$ such that the set $FS(A):=\{a_1+a_2+\dots+a_n:n\in\mathbb{N}\land a_1,a_2,\dots,a_n\in A\mbox{ are distinct}\}$ is monochromatic.
\end{te}

Let $\beta\mathbb{N}$ denote the set of ultrafilters on $\mathbb{N}$. It is common to identify, for every $n\in\mathbb{N}$, the principal ultrafilter $\{A\subseteq\mathbb{N}:n\in A\}$ with $n$, thus considering $\beta\mathbb{N}$ as an extension of $\mathbb{N}$. Observing $\mathbb{N}$ from a broader viewpoint of $\beta\mathbb{N}$ makes it easier to understand some phenomena about binary operations occurring in $\mathbb{N}$. For example, it turns out that any set belonging to a $+$-idempotent ultrafilter contains a set of the form $FS(A)$ for some infinite $A$. Since such ultrafilters exist in ZFC, this suffices to prove Hindman's finite sums theorem. The technique was invented by Galvin and Glazer, and the proof can be found in Chapter 5 of \cite{HS}. The book contains many more applications of this kind, a large number of which were found by Hindman and Strauss themselves.

Following a similar idea, to use ultrafilters to better understand the divisibility relation $\mid$ on $\mathbb{N}$, in \cite{So1} several extensions of $\mid$ were introduced. One of them turned out to have many interesting properties, some of which resemble those of divisibility on $\mathbb{N}$. Actually, extensions of relations to ultrafilters were considered in general, and in \cite{PS} two universal ways, called canonical, were isolated. One of them yields the same relation we are dealing with, although it was defined in a slightly different way. (The other one is of no use for exploring divisibility.)

This is primarily a survey paper, and therefore most of the proofs will be omitted, including instead references to corresponding papers. However, we will use the opportunity to generalize some of the theorems and include their proofs. Also, we will add a few new facts, most of them simple but connecting results that previously appeared in different papers. Sections \ref{sec3} and \ref{sec4} deal mostly with results from \cite{So3}, and Section \ref{sec5} with generalizations from \cite{So8}. Section \ref{sec6} contains results about congruence taken from \cite{So6} and \cite{DLM}, and Section \ref{sec7} is concerned with other related relations, considered in \cite{So7}.

\section{Basic definitions and notation}

What makes working with ultrafilters more accessible is the fact that the set $\beta\mathbb{N}$ of ultrafilters on $\mathbb{N}$ has a significant topological structure, when provided with base sets $\overline{A}=\{{\cal F}\in\beta\mathbb{N}:A\in{\cal F}\}$. When each $n\in\mathbb{N}$ is identified with the principal ultrafilter $\{A\subseteq\mathbb{N}:n\in A\}$, this topological space is the so-called Stone-\v Cech compactification of the discrete space on $\mathbb{N}$; the general construction is described in detail in \cite{W}. One of the features of this space is that every function $f:\mathbb{N}\rightarrow\mathbb{N}$ can be uniquely extended to a continuous $\widetilde{f}:\beta\mathbb{N}\rightarrow\beta\mathbb{N}$.

This can be used to extend also the operation $\cdot$ from $\mathbb{N}$ to $\beta\mathbb{N}$. In fact, every associative operation $*$ on $\mathbb{N}$ can be extended to an operation (also denoted by $*$) on $\beta\mathbb{N}$, making $(\beta\mathbb{N},*)$ a right-topological semigroup:
\begin{equation}
A\in{\cal F}*{\cal G}\Leftrightarrow\{n\in\mathbb{N}:n^{-1}A\in{\cal G}\}\in{\cal F},\label{extoper}
\end{equation}
where $n^{-1}A=\{m\in\mathbb{N}:n*m\in A\}$. The book \cite{HS} contains a detailed analysis of various aspects of $(\beta\mathbb{N},*)$.

Recall that the extension $\widetilde{f}$ of a function $f:\mathbb{N}\rightarrow\mathbb{N}$ was defined by setting $\widetilde{f}({\cal F})$ to be generated by direct images $f[A]=\{n\in\mathbb{N}:(\exists a\in A)f(a)=n\}$ of sets $A\in{\cal F}$. Analogously, for $A\subseteq\mathbb{N}$ let $A\uparrow:=\{n\in\mathbb{N}:(\exists a\in A)a\mid n\}$. 

\begin{de}\label{defdiv}
For ${\cal F},{\cal G}\in\beta\mathbb{N}$, ${\cal F}\widemid{\cal G}$ if and only if $(\forall A\in{\cal F})A\uparrow\in{\cal G}.$
\end{de}

For $A\subseteq\mathbb{N}$ let also $A\downarrow:=\{n\in\mathbb{N}:(\exists a\in A)n\mid a\}$. Now let ${\cal U}:=\{A\in P(\mathbb{N})\setminus\{\emptyset\}:A\uparrow=A\}$ and ${\cal V}:=\{A\in P(\mathbb{N})\setminus\{\mathbb{N}\}:A\downarrow=A\}$ be the families of all upward and all downward closed sets. It is easy to show that
$${\cal F}\widemid{\cal G}\mbox{ if and only if }{\cal F}\cap{\cal U}\subseteq{\cal G}\mbox{ if and only if }{\cal G}\cap{\cal V}\subseteq{\cal F}.$$
Since $\widemid$ is not antisymmetric, we define ${\cal F}=_\sim{\cal G}$ if and only if ${\cal F}\widemid{\cal G}$ and ${\cal G}\widemid{\cal F}$, and work with respective equivalence classes $[{\cal F}]_\sim$.

\begin{ex}
(a) Divisibility by $m\in\mathbb{N}$ is very simple to check: $m\widemid{\cal F}$ if and only if $m\mathbb{N}:=\{mn:n\in\mathbb{N}\}\in{\cal F}$. Let us call an ultrafilter $\mathbb{N}$-free if it is not divisible by any $m\in\mathbb{N}$. Several equivalent conditions for being $\mathbb{N}$-free can be found in \cite{So5}, Theorem 5.3.

(b) There are some special $=_\sim$-equivalence classes. The most prominent is the maximal class, whose existence follows from the fact that $\cal U$ has the finite intersection property. So let us denote $MAX=\{{\cal F}\in\beta\mathbb{N}:(\forall{\cal G}\in\beta\mathbb{N}){\cal G}\widemid{\cal F}\}=\{{\cal F}\in\beta\mathbb{N}:{\cal U}\subseteq{\cal F}\}$.

\end{ex}

\begin{de}
If $\{{\cal F}_i:i\in I\}$ is a family of ultrafilters and ${\cal W}$ is an ultrafilter on $I$, ${\cal G}=\lim_{i\rightarrow{\cal W}}{\cal F}_i$ is the ultrafilter defined by: $A\in{\cal G}$ if and only if $\{i\in I:A\in{\cal F}_i\}\in{\cal W}$.
\end{de}

It is shown in \cite{So4}, Example 4.2, that $(\beta\mathbb{N}/=_\sim,\widemid)$ is not well-founded.

\begin{lm}[\cite{So5}, Lemma 4.1.]\label{sublimit}
(a) Every chain $\langle[{\cal F}_i]_\sim:i\in I\rangle$ in the order $(\beta \mathbb{N}/\hspace{-1mm}=_\sim,\widemid)$ has the least upper bound $[{\cal G}_U]_\sim$ and the greatest lower bound $[{\cal G}_L]_\sim$.

(b) $\bigcup_{i\in I}({\cal F}_i\cap{\cal U})={\cal G}_U\cap{\cal U}$ and $\bigcap_{i\in I}({\cal F}_i\cap{\cal U})={\cal G}_L\cap{\cal U}$.
\end{lm}

The least upper bound $[{\cal G}_U]_\sim$ is obtained by ${\cal G}_U=\lim_{i\rightarrow{\cal W}}{\cal F}_i$ for any ultrafilter ${\cal W}$ containing all final segments of $I$, and the greatest lower bound in a similar way, using ${\cal W}$ containing all initial segments of $I$. The obtained $=_\sim$-equivalence classes do not depend on the particular $\cal W$ used.

Throughout the text, $\mathbb{N}$ will denote the set of all natural numbers (without zero) and $\omega=\mathbb{N}\cup\{0\}$. Also, $\mathbb{P}$ is the set of primes and $\mathbb{P}^{exp}:=\{p^n:p\in\mathbb{P}\land n\in\mathbb{N}\}$. For ${\cal P}\in\overline{P}$ we will denote ${\cal P}\upharpoonright\mathbb{P}=\{A\in{\cal P}:A\subseteq\mathbb{P}\}$.

\section{The nonstandard method}

To understand $\beta\mathbb{N}$ better it is often useful to parallelly consider a nonstandard extension $\zve{\mathbb{N}}$ of $\mathbb{N}$. There are several ways to define such extensions, producing models of various richness. For most of our purposes it will suffice to consider objects from a superstructure $V(X)$ containing $\mathbb{N}$, an approach introduced in \cite{RZ} and described in detail in \cite{CK}. So let $X$ be a set of atoms containing a copy of $\mathbb{N}$ (the elements of which we will nevertheless identify with natural numbers). One defines: $V_0(X)=X$, $V_{n+1}(X)=V_n(X)\cup P(V_n(X))$ for $n\in\omega$ and $V(X)=\bigcup_{n<\omega}V_n(X)$. A nonstandard extension is a superstructure $V(Y)$ such that $X\subseteq Y$, along with a star-map $*:V(X)\rightarrow V(Y)$, taking every $x\in V(X)$ into its star-counterpart $\zve x$, which is not onto and satisfies the following main tool.\\

{\it The Transfer Principle.} For every bounded formula $\varphi$ and all $a_1,a_2,\dots,a_n\in V(X)$,
\begin{equation}\label{eqtransfer}
V(X)\models\varphi(a_1,a_2,\dots,a_n)\mbox{ if and only if }V(Y)\models\varphi(\zve a_1,\zve a_2,\dots,\zve a_n).
\end{equation}
A bounded formula is a first-order formula in which all quantifiers are bounded: of the form $(\forall x\in y)$ or $(\exists x\in y)$. Note that all relations and functions on $\mathbb{N}$, on $P(\mathbb{N})$ etc.\ are actually elements of $V(X)$, so they can be considered as values of free variables in $\varphi$ on the left-hand side of (\ref{eqtransfer}), and replaced with their star-counterparts on the right-hand side.

\begin{ex}\label{extransfer}
(a) Since $(\forall x\in\mathbb{N})1\leq x$, by Transfer we have $(\forall x\in\zve{\mathbb{N}})\zve 1\hspace{1mm}\zvez\leq\hspace{1mm}x$. Analogously, Transfer implies that the elements $\zve n$ for $n\in\mathbb{N}$ constitute a proper initial part of $\zve{\mathbb{N}}$. Therefore they are identified with natural numbers, just like principal ultrafilters are. Hence $\zve{\leq}$ agrees with $\leq$ on $\mathbb{N}$.

The values of free parameters of the first formula are $1\in V_0(X)$, $\mathbb{N}\in V_1(X)$ and the relation $\leq\in V_3(X)$; in the second formula they are all replaced by their star-counterparts. It is common, however, to omit the star sign in front of the standard relation and operation simbols, such as $\leq$, to make the formulas more readable. In fact, we can safely assume the $\zve{=}$ and $\zve{\in}$ to really be the equality and set membership, see \cite{H} for an analysis of this issue.

(b) Since $(\forall x,y\in\mathbb{N})(x\mid y\Leftrightarrow(\exists k\in\mathbb{N})y=kx)$, by Transfer we get $(\forall x,y\in\zve{\mathbb{N}})(x\zvemid y\Leftrightarrow(\exists k\in\zve{\mathbb{N}})y=kx)$. In the case of divisibility relation we will disregard the rule from (a) and write $\zvemid$ for its extension. Again, this relation $\zvemid$ agrees with the usual divisibility on $\mathbb{N}$.

(c) Let $\langle p_n:n\in\mathbb{N}\rangle$ be the increasing enumeration of $\mathbb{P}$. Then its nonstandard extension $\langle p_n:n\in\zve{\mathbb{N}}\rangle$ is the increasing enumeration of $\zve{\mathbb{P}}$, the set of nonstandard primes. Recall that, for $p\in\mathbb{P}$ and $n,k\in\mathbb{N}$, $p^k\parallel n$ means that $k=\max\{l\in\mathbb{N}:p^l\mid n\}$; we say that $p^k$ is an exact divisor of $n$. Likewise, for $p\in\zve{\mathbb{P}}$, $x,k\in\zve{\mathbb{N}}$, $p^k\zvepar x$ means that $k=\max\{l\in\zve{\mathbb{N}}:p^l\zvemid x\}$.
\end{ex}

$a\in V(Y)$ is internal if it belongs to $\zve x$ for some $x\in V(X)$. To obtain a richer $\zve{\mathbb{N}}$, one can also require it to satisfy $\kappa$-saturation for some infinite cardinal $\kappa$; this means that every family $F$ of internal sets with $|F|<\kappa$ which has a finite intersection property also has a nonempty intersection. For our purposes ${\goth c}^+$-saturation will suffice, and such extensions exist in ZFC. Therefore we will make this assumption whenever needed. All these and many other aspects of nonstandard models were described in many books, of which we recommend \cite{G2} and \cite{DGL}.

A proof of the following folklore fact can be found in \cite{So4}, Theorem 2.5.

\begin{pp}\label{fund}
Let $\langle p_n:n\in\zve{\mathbb{N}}\rangle$ be the increasing enumeration of $\zve{\mathbb{P}}$ as in Example \ref{extransfer}(c).

(a) For every $z\in\zve{\mathbb{N}}$ and every internal sequence $\langle h(n):n\leq z\rangle$ there is unique $x\in\zve{\mathbb{N}}$ such that $p_n^{h(n)}\zvepar x$ for $n\leq z$ and $p_n\nzvemid x$ for $n>z$; we denote this element by $\prod_{n\leq z}p_n^{h(n)}$.

(b) Every $x\in\zve{\mathbb{N}}$ can be uniquely represented as $\prod_{n\leq z}p_n^{h(n)}$ for some $z\in\zve{\mathbb{N}}$ and some internal sequence $\langle h(n):n\leq z\rangle$ such that $h(z)>0$.
\end{pp}

The connection between nonstandard natural numbers and ultrafilters on $\mathbb{N}$ was described in \cite{Lux}. Here we change some of our notation from previous papers to one with a more model-theoretic flavour, as suggested in \cite{DLM}. For every $x\in\zve{\mathbb{N}}$, the family $\{A\subseteq\mathbb{N}:x\in\zve A\}$ is an ultrafilter ${\cal F}$; we will write $x\models{\cal F}$ or ${\cal F}=tp(x/\zve{\mathbb{N}})$ and say that $x$ is a generator of $\cal F$. In a ${\goth c}^+$-saturated nonstandard extension the set of generators $\mu({\cal F}):=\{x\in\zve{\mathbb{N}}:x\models{\cal F}\}$ is nonempty for every ${\cal F}\in\beta\mathbb{N}$. Also, if $F=\{{\cal F}_i:i\in I\}$ is a family of ultrafilters, let $\mu(F)=\bigcup_{i\in I}\mu({\cal F}_i)$.

As noticed in \cite{NR}, the mapping $x\mapsto tp(x/\zve{\mathbb{N}})$ is interchangeable with the extensions of a function $f:\mathbb{N}\rightarrow\mathbb{N}$ in the following sense: for every $x\in\zve{\mathbb{N}}$,
\begin{equation}\label{eqv}
tp(\zve f(x)/\zve{\mathbb{N}})=\widetilde{f}(tp(x/\zve{\mathbb{N}})).
\end{equation}

If $x\models{\cal F}$ and $y\models{\cal G}$, then the pair $(x,y)$ is a generator of an ultrafilter on $\mathbb{N}\times\mathbb{N}$. However, this ultrafilter does not need to be the tensor product
$${\cal F}\otimes{\cal G}:=\{A\subseteq\mathbb{N}\times\mathbb{N}:\{m\in\mathbb{N}:\{n\in\mathbb{N}:(m,n)\in A\}\in{\cal G}\}\in{\cal F}\}.$$

\begin{de}
If $x\models{\cal F}$, $y\models{\cal G}$ and $(x,y)\models{\cal F}\otimes{\cal G}$, then $(x,y)$ is called a tensor pair.
\end{de}

Tensor pairs were used informally by Puritz, and subsequently named and studied in detail by Di Nasso in \cite{DN2}. In \cite{L4} Luperi Baglini generalized the notion to tensor $k$-tuples, and gave many equivalent conditions.

Since any bijection $\varphi:\mathbb{N}\times\mathbb{N}\rightarrow\mathbb{N}$ yields a homeomorphism $\widetilde{\varphi}:\beta(\mathbb{N}\times\mathbb{N})\rightarrow\beta\mathbb{N}$, we can apply the analogue of (\ref{eqv}) to the multiplication or the addition function $f:\mathbb{N}\times\mathbb{N}\rightarrow\mathbb{N}$; in this way we get that,  if $(x,y)$ is a tensor pair, then $tp(x\cdot y/\zve{\mathbb{N}})={\cal F}\cdot{\cal G}$ and $tp(x+y/\zve{\mathbb{N}})={\cal F}+{\cal G}$.

\begin{te}[\cite{DN2}, Theorem 11.5.12]\label{tensor2}
For every $x\in\zve{\mathbb{N}}\setminus\mathbb{N}$ and every ${\cal G}\in\beta\mathbb{N}\setminus\mathbb{N}$, there are $y,y'\in\mu({\cal G})$ such that $(x,y)$ and $(y',x)$ are tensor pairs.
\end{te}

\begin{te}[\cite{P}, Theorem 3.4]\label{tensor}
$(x,y)\in\zve{\mathbb{N}}\times\zve{\mathbb{N}}$ is a tensor pair if and only if, for every $f:\mathbb{N}\rightarrow\mathbb{N}$, either $\zve f(y)\in\mathbb{N}$ or $\zve f(y)>x$.
\end{te}

By Lemma \ref{sublimit}, every $\widemid$-chain has the smallest upper bound and the greatest lower bound. As the next result shows, chains in this order are reflected in chains in $(\zve{\mathbb{N}},\zvemid)$.

\begin{te}[\cite{So5}, Lemma 4.4]\label{chains}
Let $V(Y)$ be ${\goth c}^+$-saturated. For any well-ordered strictly increasing $\widemid$-chain $\langle{\cal G}_\alpha:\alpha<\beta\rangle$ there is a $\zvemid$-chain $\langle x_\alpha:\alpha<\beta\rangle$ such that $x_\alpha\models{\cal G}_\alpha$ for $\alpha<\beta$.
\end{te}

There are several arguments suggesting that $\widemid$ might be "the right" way to extend the divisibility relation to $\beta\mathbb{N}$. One is the fact that it is the same as the canonical extension from \cite{PS}. Another is the following theorem, saying that the $\widemid$-divisibility of ultrafilters reflects the $\zvemid$-divisibility on $\zve{\mathbb{N}}$. 

\begin{te}[\cite{So5}, Theorem 3.4]\label{ekviv}
For every ${\goth c}^+$-saturated extension and every two ultrafilters ${\cal F},{\cal G}\in\beta\mathbb{N}$ the following conditions are equivalent:

(i) ${\cal F}\widemid{\cal G}$;

(ii) there are $x\models{\cal F}$ and $y\models{\cal G}$ such that $x\zvemid y$;

(iii) for every $x\models{\cal F}$ there is $y\models{\cal G}$ such that $x\zvemid y$;

(iv) for every $y\models{\cal G}$ there is $x\models{\cal F}$ such that $x\zvemid y$.
\end{te}

Another argument is that a significant part of $(\beta\mathbb{N}/=_\sim,\widemid)$ resembles $(\mathbb{N},\mid)$ in many ways, and we will consider this resemblance in the next section.

\section{The lower part $L$}\label{sec3}

For $n\in\omega$ denote $L_n=\{a_1a_2\dots a_n:a_1,a_2,\dots,a_n\in\mathbb{P}\}$. Thus, $\overline{L_0}=L_0=\{1\}$ and $L_1=\mathbb{P}$. We will see that $L:=\bigcup_{n<\omega}L_n$, "the lower half" of $(\beta\mathbb{N}/=_\sim,\widemid)$, is divided into levels $\overline{L_n}$ resembling the levels $L_n$ in the order $(\mathbb{N},\mid)$.

\begin{de}
An ultrafilter ${\cal P}\in\beta\mathbb{N}\setminus\{1\}$ is prime if it is $\widemid$-divisible only by $1$ and itself.
\end{de}

\begin{te}[\cite{So3}, Theorem 2.3]\label{widemin}
${\cal P}\in\beta\mathbb{N}$ is prime if and only if ${\cal P}\in\overline{\mathbb{P}}$.
\end{te}

\begin{te}[\cite{So5}, Theorem 2.2]
If ${\cal P}\in\overline{\mathbb{P}}$ and ${\cal F},{\cal G}\in\beta\mathbb{N}$, then ${\cal P}\widemid{\cal F}\cdot{\cal G}$ implies ${\cal P}\widemid{\cal F}$ or ${\cal P}\widemid{\cal G}$.
\end{te}

Our goal now is to describe how the place of an ultrafilter in the $\widemid$-hierarchy depends on its basic ingredients. For these ingredients we can not, however, take only prime ultrafilters, but also their powers.

If $f_n:\mathbb{N}\rightarrow\mathbb{N}$ is defined by $f_n(m)=m^n$, for ${\cal F}\in\beta\mathbb{N}$ let us denote $\widetilde{f_n}({\cal F})$ by ${\cal F}^n$. Unlike divisibility in $\mathbb{N}$, there are ultrafilters divisible by only one prime $\cal P$ which are different from powers of $\cal P$. For this reason we need some more definitions.

\begin{de}\label{patternL}
For $A,B,A_i\subseteq\mathbb{N}$ and $n\in\mathbb{N}$, let:
\begin{eqnarray*}
A^n &=& \{a^n:a\in A\}\\
A_1A_2\dots A_n &=& \{a_1a_2\dots a_n:a_i\in A_i\mbox{ for }1\leq i\leq n\land gcd(a_i,a_j)=1\mbox{ for }i\neq j\}\\
A^{(n)} &=& \underbrace{A\cdot A\cdot\dots\cdot A}_n.
\end{eqnarray*}
We call ultrafilters of the form ${\cal P}^k$ for some ${\cal P}\in\overline{\mathbb{P}}$ and $k\in\mathbb{N}$ {\it basic}. Let ${\cal B}$ be the set of all basic ultrafilters, and let ${\cal A}$ be the set of all functions $\alpha:{\cal B}\rightarrow\omega$ with finite support (i.e.\ such that ${\rm supp}(\alpha):=\{{\cal P}^k\in{\cal B}:\alpha({\cal P}^k)\neq 0\}$ is finite); we call functions $\alpha\in{\cal A}$ patterns.
\end{de}

\begin{ex}\label{transprod}
More easy applications of the Transfer principle show that, for $A,B\in P(\mathbb{N})$, $k\in\mathbb{N}$ and $f:\mathbb{N}\rightarrow\mathbb{N}$: (a) $\zve(A^k)=(\zve A)^k$; (b) $\zve(AB)=\zve A\zve B$ and (c) $\zve(f[A])=\zve f[\zve A]$. We will use these properties extensively in the rest of the paper.
\end{ex}

We will abuse notation and write $\alpha=\{({\cal P}_1^{k_1},m_1),({\cal P}_2^{k_2},m_2),\dots,({\cal P}_n^{k_n},m_n)\}$ if $\alpha({\cal P}^{k})=0$ for ${\cal P}^{k}\in{\cal B}\setminus\{{\cal P}_i^{k_i}:1\leq i\leq n\}$ (allowing also some of the $m_i$ to be zeros).

Now we want to adjoin to each ${\cal F}\in\beta\mathbb{N}$ a pattern and identify the upward closed sets belonging to $\cal F$ which are determined by its pattern. They will be exactly the upward-closures of the sets described in the following definition.

\begin{de}\label{defFalpha}
Let $\alpha=\{({\cal P}_1^{k_1},m_1),({\cal P}_2^{k_2},m_2),\dots,({\cal P}_n^{k_n},m_n)\}\in{\cal A}$ (${\cal P}_i\in\overline{\mathbb{P}}$). With $F_\alpha$ we denote the family of all sets
\begin{eqnarray*}\label{eqhigher}
&& (A_1^{k_1})^{(m_1)}(A_2^{k_2})^{(m_2)}\dots(A_n^{k_n})^{(m_n)}\\
&=& \left\{\prod_{i=1}^n\prod_{j=1}^{m_i}p_{i,j}^{k_i} : p_{i,j}\in A_i\mbox{ for all }i,j\;\land\mbox{ all }p_{i,j}\mbox{ are distinct}\right\}
\end{eqnarray*}
such that: (i) $A_i\in{\cal P}_i\upharpoonright\mathbb{P}$, (ii) $A_i=A_j$ if ${\cal P}_i={\cal P}_j$ and $A_i\cap A_j=\emptyset$ otherwise.
\end{de}

\begin{ex}
(a) For $\alpha=\{({\cal P}^k,1)\}$, the only ultrafilter containing $F_\alpha$ is ${\cal P}^k$ itself.

(b) For $\alpha=\{({\cal P},1),({\cal Q},1)\}$, the ultrafilters containing $F_\alpha$ are exactly those divisible only by $1,{\cal P},{\cal Q}$ and themselves. One such ultrafilter is ${\cal P}\cdot{\cal Q}$, and another (usually different) is ${\cal Q}\cdot{\cal P}$. This will be generalized in Theorem \ref{products}.

(c) If $\alpha=\{({\cal P},2)\}$, then $F_\alpha=\{A^{(2)}:A\in{\cal P}\upharpoonright\mathbb{P}\}$. Ultrafilters containing $F_\alpha$ are divisible only by $1,{\cal P}$ and themselves. An example is ${\cal P}\cdot{\cal P}$, which is distinct from ${\cal P}^2$.

(d) Finally, if $\alpha=\{({\cal P}^2,1),({\cal P},2),({\cal Q},1)\}$, the corresponding ultrafilters contain sets of the form $A^2A^{(2)}B$ for disjoint $A\in{\cal P}$, $B\in{\cal Q}$. We will see in what follows that such ultrafilters belong to the fifth level of the $\widemid$-hierarchy. Their only divisors on the first level are $\cal P$ and $\cal Q$. Their divisors on the second level are ${\cal P}^2$, some ultrafilters with pattern $\{({\cal P},2)\}$ as in (c), and some ultrafilters with pattern $\{({\cal P},1),({\cal Q},1)\}$. They also have divisors on levels 3 and 4, described in a similar manner.
\end{ex}

\begin{de}
If $\alpha=\{({\cal P}_1^{k_1},m_1),({\cal P}_2^{k_2},m_2),\dots,({\cal P}_n^{k_n},m_n)\}\in{\cal A}$ (${\cal P}_i\in\overline{\mathbb{P}}$), we denote $\sigma(\alpha)=\sum_{i=1}^nk_im_i$.
\end{de}

In \cite{So3} a proof of the special case $n=4$ of the following theorem was given, containing all the essential ingredients. We now give a complete proof.

\begin{te}[\cite{So3}, Theorem 5.5]\label{nivoi}
The set $\overline{L_n}$ (for $n\in\mathbb{N}$) consists precisely of ultrafilters containing $F_\alpha$ for some $\alpha\in{\cal A}$ such that $\sigma(\alpha)=n$.
\end{te}

\dokaz If $\sigma(\alpha)=n$, from Definition \ref{defFalpha} it is clear that any set $(A_1^{k_1})^{(m_1)}(A_2^{k_2})^{(m_2)}$ $\dots(A_n^{k_n})^{(m_n)}\in F_\alpha$ is a subset of $L_n$, so any ultrafilter containing sets from $F_\alpha$ also contains $L_n$.

Now let $X_n$ be the family of all functions $\varphi:\mathbb{N}\rightarrow\omega$ such that $\sum_{k\in\mathbb{N}}k\varphi(k)=n$. Clearly, for each $\varphi\in X_n$ the support ${\rm supp}\varphi:=\{k\in\mathbb{N}:\varphi(k)\neq 0\}$ is finite. Thus,
$$L_n=\bigcup_{\varphi\in X_n}\prod_{k\in{\rm supp}\varphi}(\mathbb{P}^k)^{(\varphi(k))}.$$
Since every $X_n$ is finite, each ${\cal F}\in\overline{L_n}$ contains $Q:=\prod_{k\in{\rm supp}\varphi}(\mathbb{P}^k)^{(\varphi(k))}$ for some $\varphi\in X_n$. We define functions $f_{k,j}:Q\rightarrow\mathbb{P}$ for $k\in{\rm supp}\varphi$ and $1\leq j\leq\varphi(k)$ as follows: if $x=\prod_{k\in{\rm supp}\varphi}(p_{k,1}^kp_{k,2}^k\dots p_{k,\varphi(k)}^k)\in Q$, where $p_{k,1}<p_{k,2}<\dots<p_{k,\varphi(k)}$, let $f_{k,j}(x)=p_{k,j}$. Then $\widetilde{f_{k,j}}({\cal F})$ is an ultrafilter in $\overline{\mathbb{P}}$; denote it by ${\cal P}_{k,j}$. For ${\cal Q}^k\in{\cal B}$ define $\alpha({\cal Q}^k)=|\{j\leq\varphi(k):{\cal P}_{k,j}={\cal Q}\}|$.

To prove that ${\cal F}\supseteq F_\alpha$, let $B:=(A_1^{k_1})^{(m_1)}(A_2^{k_2})^{(m_2)}\dots(A_n^{k_n})^{(m_n)}\in F_\alpha$ be as in Definition \ref{defFalpha}. For each $l=1,2,\dots,n$ there are $m_l$ indices $j_1<j_2<\dots<j_{m_l}$ such that $A_l\in{\cal P}_{k_l,j_1}={\cal P}_{k_l,j_2}=\dots={\cal P}_{k_l,j_{m_l}}$. Then, for $r=1,2,\dots,m_l$, $\widetilde{f_{{k_l},{j_r}}}({\cal F})={\cal P}_{k_l,j_r}$ so ${f_{{k_l},{j_r}}}^{-1}[A_l]\in{\cal F}$. The intersection of these sets (for all $l=1,2,\dots,n$ and all $r=1,2,\dots,m_l$) is a subset of $B$ (more precisely, this intersection consists of those elements $\prod_{k\in{\rm supp}\varphi}\prod_{j=1}^{\varphi(k)}p_{k,j}^k\in B$ for which, if $p_{k,1}<p_{k,2}<\dots<p_{k,\varphi(k)}$, then exactly $p_{k,j_1},p_{k,j_2},\dots,p_{k,j_{m_l}}$ belong to $A_l$), and hence $B\in{\cal F}$.\kraj

For distinct patterns $\alpha$ and $\beta$, there are disjoint sets $X\in F_\alpha$ and $Y\in F_\beta$. Therefore, each ultrafilter can contain at most one family of the form $F_\alpha$; by Theorem \ref{nivoi} every ${\cal F}\in L$ really does contain one. Denote the corresponding pattern by $\alpha_{\cal F}$.

Note also that, for $p\in\mathbb{P}$: (1) $\alpha_{\cal F}(p^k)\leq 1$ for all $k\in\mathbb{N}$ and (2) $\alpha_{\cal F}(p^k)=1$ for at most one value of $k$. Namely, if ${\cal F}=p^kp^l{\cal G}$, then $\cal F$ would actually be divisible by $p^{k+l}$.

Now let us prove a nonstandard characterization of ultrafilters containing $F_\alpha$, which will later be used to generalize the notion of a pattern from $L$ to the whole $\beta\mathbb{N}$.

\begin{te}\label{patternsL}
Let ${\cal F}\in L$ and $\alpha_{\cal F}=\{({\cal P}_1^{k_1},m_1),({\cal P}_2^{k_2},m_2),\dots,({\cal P}_n^{k_n},m_n)\}$. The generators of $\cal F$ are precisely $x\in\zve{\mathbb{N}}$ of the form
$$x=\prod_{i=1}^n\prod_{j=1}^{m_i}p_{i,j}^{k_i},$$
where $p_{i,j}\models{\cal P}_i$ for all $i,j$ and all $p_{i,j}$ are distinct.
\end{te}

\dokaz Take any set $B:=(A_1^{k_1})^{(m_1)}(A_2^{k_2})^{(m_2)}\dots(A_n^{k_n})^{(m_n)}\in F_{\alpha_{\cal F}}$ as in Definition \ref{defFalpha}. Every generator $x$ of $\cal F$ belongs to $\zve B=(\zve A_1^{k_1})^{(m_1)}(\zve A_2^{k_2})^{(m_2)}\dots(\zve A_n^{k_n})^{(m_n)}$ (see Example \ref{transprod}), so it must be of the form $x=\prod_{i=1}^n\prod_{j=1}^{m_i}p_{i,j}^{k_i},$
where $p_{i,j}\in\zve A_i$ for all $i,j$ and all $p_{i,j}$ are distinct.

Since this is true for any set $B$ as above, and the sets $A_i$ corresponding to distinct prime ultrafilters must be disjoint, it follows that $p_{i,j}\in\zve A_i$ for all $A_i\in{\cal P}_i$, so $p_{i,j}$ is a generator of ${\cal P}_i$.

In the other direction, if $x$ is as in the formulation, then $p_{i,j}\models{\cal P}_i$ implies $p_{i,j}\in\zve A_i$, so $x\in\zve B$ for all $B\in F_{\alpha_{\cal F}}$.\kraj

Define an operation on $\cal A$ as follows. For $p\in\mathbb{P}$ denote temporarily by ${\rm exp}_p\alpha$ the (unique) $k\in\mathbb{N}$ such that $\alpha(p^k)=1$ if it exists, and ${\rm exp}_p\alpha=0$ otherwise. For $\alpha,\beta\in{\cal A}$ let ${\rm exp}_p(\alpha\oplus\beta)={\rm exp}_p\alpha+{\rm exp}_p\beta$ if $p\in\mathbb{P}$, and
$$(\alpha\oplus\beta)({\cal P}^k)=\alpha({\cal P}^k)+\beta({\cal P}^k)\mbox{ if }{\cal P}\in\overline{\mathbb{P}}\setminus\mathbb{P}\mbox{ and }k\in\mathbb{N}.$$
The following is a slight generalization of \cite{So3}, Theorem 5.14, with a simplified proof.

\begin{lm}\label{propattern}
For any ${\cal F},{\cal G}\in\beta\mathbb{N}$: 
\begin{equation}\label{eqpropattern}
\alpha_{{\cal F}\cdot{\cal G}}=\alpha_{\cal F}\oplus\alpha_{\cal G}.
\end{equation}
\end{lm}

\dokaz Let ${\cal F}=m{\cal F}'$ and ${\cal G}=n{\cal G}'$ for some $\mathbb{N}$-free ultrafilters ${\cal F}'$ and ${\cal G}'$. Clearly, the exponent of any $p\in\mathbb{P}$ in ${\cal F}\cdot{\cal G}$ is obtained by adding exponents of $p$ in $m$ and $n$. It remains to check the result for ${\cal F}'\cdot{\cal G}'$. Let $x\models{\cal F}'$ and $y\models{\cal G}'$ be such that $(x,y)$ is a tensor pair (see Theorem \ref{tensor2}). Applying Theorem \ref{tensor} to functions of the form $f_k(p_1^{a_1}p_2^{a_2}\dots p_s^{a_s})=p_k$ for $p_1<p_2<\dots<p_s$ (returning the $k$-th smallest prime if it exists, and 1 otherwise), we conclude that every prime factor of $y$ is greater than any prime factor of $x$. Since $x\cdot y\models{\cal F}\cdot{\cal G}$, (\ref{eqpropattern}) follows by Theorem \ref{patternsL}.\kraj

\begin{te}\label{products}
Let $\alpha=\{({\cal P}_1^{k_1},m_1),({\cal P}_2^{k_2},m_2),\dots,({\cal P}_n^{k_n},m_n)\}\in{\cal A}$ and ${\cal F}$ is a product of 
$$\underbrace{{\cal P}_1^{k_1},{\cal P}_1^{k_1},\dots,{\cal P}_1^{k_1}}_{m_1},\dots,\underbrace{{\cal P}_n^{k_n},{\cal P}_n^{k_n},\dots,{\cal P}_n^{k_n}}_{m_n}$$
in any order. Then $\alpha_{\cal F}=\alpha$, so ${\cal F}\supseteq F_\alpha$.
\end{te}

\dokaz For the special case $\alpha=\{({\cal P}^{k},1)\}$ the theorem is trivially true. The general case then follows using Lemma \ref{propattern}: if $\cal F$ is a product of ultrafilters, then the pattern of $\cal F$ is obtained by applying the operation $\oplus$ to patterns of factors.\kraj



\section{Levels, patterns and $=_\sim$-classes in $L$}\label{sec4}

To begin with, the ultrafilters in $L$ are completely determined by their $=_\sim$-equivalence class.

\begin{te}[\cite{So3}, Corollary 5.10]\label{singletons}
$|[{\cal F}]_\sim|=1$ for every $n\in\mathbb{N}$ and all ultrafilters ${\cal F}\in L$.
\end{te}

Let us now turn to the question: how many ultrafilters with the same pattern are there? Since, for every $\alpha\in{\cal A}$, the family $\{{\cal F}\in\beta\mathbb{N}:F_\alpha\subseteq{\cal F}\}$ is closed in $\beta\mathbb{N}$, it must be either finite or of cardinality $2^{\goth c}$.

Let ${\cal P}\in\overline{\mathbb{P}}\setminus\mathbb{P}$ and $\alpha=\{({\cal P},2)\}$. In \cite{So3}, Theorem 3.6 we proved that there is a unique ultrafilter ${\cal F}\supseteq F_\alpha$ if and only if $\cal P$ is Ramsey. Let us generalize this a little bit.

\begin{de}
${\cal F}\in\beta\mathbb{N}\setminus\mathbb{N}$ is $n$-Ramsey if $n$ is the smallest natural number such that $\omega\rightarrow[{\cal F}]^2_{n+1,\leq n}$: for every coloring $c:[\mathbb{N}]^2\rightarrow n+1$ of pairs of elements of $\mathbb{N}$ there is $A\in{\cal F}$ such that $|c[[A]^2]|\leq n$. 

${\cal F}$ is Ramsey if it is 1-Ramsey. ${\cal F}$ is weakly Ramsey if it is $n$-Ramsey for some $n>1$.
\end{de}



\begin{te}\label{ramsey}
Let ${\cal P}\in\overline{\mathbb{P}}\setminus\mathbb{P}$, $\alpha=\{({\cal P},2)\}$ and $n\in\mathbb{N}$. There are exactly $n$ ultrafilters ${\cal F}\supseteq F_\alpha$ if and only if $\cal P$ is $n$-Ramsey.
\end{te}

\dokaz Assume $\cal P$ is $n$-Ramsey. Then there is a coloring $c:[\mathbb{P}]^2\rightarrow\{0,1,\dots,n-1\}$ such that, for every $A\in{\cal P}\upharpoonright\mathbb{P}$, $[A]^2$ intersects every $c^{-1}[\{i\}]$. Hence every $A^{(2)}$ intersects every $S_i:=\{ab:c(\{a,b\})=i\}$. This means that, for every $i<n$, $F_\alpha\cup\{S_i\}$ has the finite intersection property. Thus there are ultrafilters ${\cal F}_i$ containing $F_\alpha\cup\{S_i\}$ respectively, and thus being distinct from each other.

Now assume that there are $n+1$ distinct ultrafilters ${\cal F}_0,{\cal F}_1$, $\dots,{\cal F}_{n}\supseteq F_\alpha$. Then there is a partition $\{S_0,S_1,\dots,S_{n}\}$ of $L_2$ such that $S_i\in{\cal F}_i$ for every $i\leq n$. Define a coloring of $[\mathbb{P}]^2$ as follows: let $c(\{a_1,a_2\})=i$ if $a_1a_2\in S_i$. Since $\cal P$ is $n$-Ramsey, there is a set $A\in{\cal P}$ such that some $i\notin c[[A]^2]$. This means that $A^{(2)}\cap S_i=\emptyset$, so $\alpha$ can not be the pattern of ${\cal F}_i$.\kraj

The next two results show that the case of $2^{\goth c}$-many ultrafilters of the same pattern is also possible.

\begin{te}[\cite{So3}, Theorem 3.13]
Assume CH. Then there is ${\cal P}\in\overline{\mathbb{P}}\setminus\mathbb{P}$ such that for every $n\geq 2$ there are $2^{\goth c}$ ultrafilters ${\cal F}$ such that $\alpha_{\cal F}=\{({\cal P},n)\}$.
\end{te}

\begin{te}[\cite{So3}, Theorem 4.5]
For every ${\cal P}\in\overline{\mathbb{P}}\setminus\mathbb{P}$ there is an ultrafilter ${\cal Q}\in\overline{\mathbb{P}}\setminus\mathbb{P}$ such that there are $2^{\goth c}$ ultrafilters ${\cal F}$ with $\alpha_{\cal F}=\{({\cal P},1),({\cal Q},1)\}$.
\end{te}

\begin{te}[\cite{So3}, Theorem 4.4]
Let ${\cal P},{\cal Q}\in\overline{\mathbb{P}}\setminus\mathbb{P}$. If there is unique ${\cal F}$ such that $\alpha_{\cal F}=\{({\cal P},1),({\cal Q},1)\}$ then both ${\cal P}$ and ${\cal Q}$ are P-points.
\end{te}

Limits of ultrafilters were first considered in the context of divisibility in order to help the understanding of ultrafilters outside $L$. We conclude this section with several examples of behavior of elements of $L$ under limits. We also include a version of sum of ultrafilters:

\begin{de}
For ${\cal F},{\cal G}_n,{\cal K}\in\beta\mathbb{N}$, we write ${\cal K}=\sum_{n\rightarrow{\cal F}}{\cal G}_n$ if: $A\in{\cal K}$ whenever $\{n\in\mathbb{N}:A/n\in{\cal G}_n\}\in{\cal F}$.
\end{de}

We remark that in \cite{F} by the sum of ultrafilters Frolik meant what is today called the limit of ultrafilters, and that often the sum is nowadays defined in a similar manner as ours, but producing an ultrafilter on $\mathbb{N}\times\mathbb{N}$ as result. Our version can, however, be obtained from limits as follows: $\sum_{n\rightarrow{\cal F}}{\cal G}_n=\{A\subseteq\mathbb{N}:\{n\in\mathbb{N}:A\in n{\cal G}_n\}\in{\cal F}\}=\lim_{n\rightarrow{\cal F}}(n{\cal G}_n)$. The product ${\cal F}\cdot{\cal G}$ from (\ref{extoper}) can also be written as this kind of sum: if ${\cal G}_n={\cal G}$ for all $n\in\mathbb{N}$, then $\sum_{n\rightarrow{\cal F}}{\cal G}=\lim_{n\rightarrow{\cal F}}(n{\cal G})=(\lim_{n\rightarrow{\cal F}}n){\cal G}={\cal F}\cdot{\cal G}$, since multiplication by ${\cal G}$ from the right is continuous.

\begin{ex}
Let $\{{\cal G}_n:n\in\mathbb{N}\}\subseteq\beta\mathbb{N}$, ${\cal H}=\lim_{n\rightarrow{\cal F}}{\cal G}_n$ and ${\cal K}=\sum_{n\rightarrow{\cal F}}{\cal G}_n$. Then:

(a) ${\cal F}\widemid{\cal K}$: for every $A\in{\cal F}\cap{\cal U}$ and every $n\in A$, $A/n=\mathbb{N}\in{\cal G}_n$, so $\{n\in\mathbb{N}:A/n\in{\cal G}_n\}\supseteq A\in{\cal F}$ and $A\in{\cal K}$. 

${\cal H}\widemid{\cal K}$: if $A\in{\cal H}\cap{\cal U}$, then $A\subseteq A/n$, so $\{n\in\mathbb{N}:A\in{\cal G}_n\}\in{\cal F}$ implies $\{n\in\mathbb{N}:A/n\in{\cal G}_n\}\in{\cal F}$, i.e.\ $A\in{\cal K}$.

(b) Let ${\cal F}\in\overline{L_m}$ and ${\cal G}_n\in\overline{L_k}$ for all $n\in\mathbb{N}$. $L_k\in{\cal G}_n$ for all $n\in\mathbb{N}$ implies $\{n\in\mathbb{N}:L_k\in{\cal G}_n\}=\mathbb{N}\in{\cal F}$, so $L_k\in{\cal H}$ and ${\cal H}\in\overline{L_k}$ as well.

${\cal K}\in\overline{L_{m+k}}$: for every $n\in L_m$, $L_{m+k}/n=L_k$. Hence, $\{n\in\mathbb{N}:L_{m+k}/n\in{\cal G}_n\}\supseteq L_m\in{\cal F}$, so $L_{m+k}\in{\cal K}$.
\end{ex}

\section{Generalized patterns}\label{sec5}

In this section we use the characterization of generators obtained in Theorem \ref{patternsL} to generalize patterns of ultrafilters. First we need to generalize basic classes to include infinite "powers" of prime ultrafilters.

\begin{de}
Let ${\cal P}\in\overline{\mathbb{P}}\setminus\mathbb{P}$. The relation $\approx_{\cal P}$ on $\mu({\cal P})\times\zve{\mathbb{N}}$ is defined as follows:
$$(p,x)\approx_{\cal P}(q,y)\mbox{ if and only if }tp(p^x/\zve{\mathbb{N}})=_\sim tp(q^y/\zve{\mathbb{N}}).$$
$\approx_{\cal P}$ is an equivalence relation, so let
$${\cal E}_{\cal P}=\{[(p,x)]_{\approx_{\cal P}}:(p,x)\in\mu({\cal P})\times\zve{\mathbb{N}}\}$$
be the set of its equivalence classes. For any $p\in\mu({\cal P})$ and $u\in{\cal E}_{\cal P}$ let $u_p:=\{x\in\zve{\mathbb{N}}:(p,x)\in u\}$.

Families of ultrafilters of the form ${\cal P}^u:=\{tp(p^x/\zve{\mathbb{N}}):(p,x)\in u\}$ for some ${\cal P}\in\overline{\mathbb{P}}$ and $u\in{\cal E}_{\cal P}$ will be called basic.

By ${\cal B}$ we denote the set of all basic classes.
\end{de}

Of course, this $\cal B$ expands the family $\cal B$ from Definition \ref{patternL}. For $k\in\mathbb{N}$ we identify the class $[(p,x)]_{\approx_{\cal P}}$ with $k$. In general, if ${\cal P}^u=\{{\cal P}^k\}$ is a singleton, we identify ${\cal P}^u$ with ${\cal P}^k$ and write $\alpha({\cal P}^k)$ instead of $\alpha({\cal P}^u)$.

\begin{lm}[\cite{So8}, Lemma 2.5]\label{approxosob}
Let ${\cal P}\in\overline{\mathbb{P}}\setminus\mathbb{P}$ and $u\in{\cal E}_{\cal P}$.

(a) All elements of $\mu({\cal P}^u)$ are of the form $p^x$ for some $(p,x)\in u$.

(b) For every $p\in\mu({\cal P})$ the set $u_p$ is nonempty and convex: if $x,y\in u_p$ and $x<z<y$, then $z\in u_p$ as well.

(c) Each $u_p$ is either a singleton or a union of galaxies.
\end{lm}

It is easy to see that if, for some $p\models{\cal P}$, $u_p$ is a singleton, then ${\cal P}^u$ is a singleton, so $u_q$ is also a singleton for all $q\models{\cal P}$. In this case we say that the basic class ${\cal P}^u$ is {\it of the first kind}, and if $u_p$ is a union of galaxies for all $p\models{\cal P}$, we will say that ${\cal P}^u$ is {\it of the second kind}. It is shown in \cite{So8} that there are basic classes of both kinds outside $L$. Note that classes of the first kind give us more examples of singleton $=_\sim$-equivalence classes (recall that in $L$ all classes are singletons by Theorem \ref{singletons}).

\begin{de}
On ${\cal E}_{\cal P}$ we define the relation:
\begin{eqnarray*}
u\prec_{\cal P}v &\mbox{ if and only if } & u\neq v\mbox{ and for some } p\in\mu({\cal P})\mbox{ and some }x,y\in\zve{\mathbb{N}}\\
&& \mbox{ holds }(p,x)\in u,(p,y)\in v\mbox{ and }x<y.
\end{eqnarray*}
We write $u\preceq_{\cal P}v$ if $u\prec_{\cal P}v$ or $u=v$.
\end{de}

As shown in \cite{So4}, Example 4.5, for $p\in\mathbb{P}$ the order $({\cal E}_p,\prec_p)$ is isomorphic to $\omega+1$: there is only one "infinite" basic class. Things are different for ${\cal P}\in\overline{\mathbb{P}}\setminus\mathbb{P}$.

\begin{lm}[\cite{So8}, Lemma 2.7]\label{supinf}
For every ${\cal P}\in\overline{\mathbb{P}}\setminus\mathbb{P}$:

(a) $\prec_{\cal P}$ is a strict linear order.

(b) Every increasing sequence in $({\cal E}_{\cal P},\prec_{\cal P})$ has a supremum and every decreasing sequence has an infimum.

(c) The order $({\cal E}_{\cal P},\prec_{\cal P})$ contains a copy of $(\mathbb{R},<)$.
\end{lm}

In particular, every ${\cal E}_{\cal P}$ has a maximum, which we denote by ${\cal P}^{max}$.

In \cite{So8} we asked whether the orders $({\cal E}_{\cal P},\prec_{\cal P})$ are isomorphic for all ${\cal P}\in\overline{\mathbb{P}}\setminus\mathbb{P}$. In \cite{DLMS} it will be shown that the answer to this question is independent of ZFC.

Let us now apply these generalized primes to answer Question 4.3 from \cite{So5}.

\begin{te}
There is an $=_\sim$-equivalence class $[{\cal F}]_\sim\in\beta\mathbb{N}/=_\sim$ which can not be represented as a limit of a $\widemid$-increasing chain.
\end{te}

\dokaz Let ${\cal P}^u=\{{\cal F}\}$ be a basic class of the first kind such that ${\cal F}\notin L$. Assume there is a $\widemid$-increasing chain $\langle[{\cal G}_\alpha]_\sim:\alpha<\beta\rangle$ whose limit is $[\cal F]_\sim$. (If the chain were not well-ordered, we could extract a well-ordered cofinal subchain and work with it instead.) This means, by Lemma \ref{sublimit}, that ${\cal F}=\lim_{\alpha\rightarrow{\cal W}}{\cal G}_\alpha$ for some ultrafilter $\cal W$ on $\beta$. By Theorem \ref{chains}, if we work with a ${\goth c}^+$-saturated nonstandard extension, we can find a $\zvemid$-chain $\langle x_\alpha:\alpha\leq\beta\rangle$ such that $x_\alpha\models{\cal G}_\alpha$ for $\alpha<\beta$ and $x_\beta\models{\cal F}$. By Lemma \ref{approxosob}, $x_\beta=p^{z_\beta}$ for some $p\in\zve{\mathbb{P}}$ and some $z_\beta\in\zve{\mathbb{N}}$. Hence there is an increasing sequence $\langle z_\alpha:\alpha\leq\beta\rangle$ such that $x_\alpha=p^{z_\alpha}$ for all $\alpha<\beta$. However, since ${\cal P}^u$ is of the first kind, $p^{z_\beta}$ is the only power of $p$ which is a generator of $\cal F$, so $p^{z_\beta-1}$ is also an upper bound of $\langle x_\alpha:\alpha<\beta\rangle$. Thus, $[tp(p^{z_\beta-1}/\zve{\mathbb{N}})]_\sim$ is a upper bound of $\langle[{\cal G}_\alpha]_\sim:\alpha<\beta\rangle$ strictly smaller than $[\cal F]_\sim$, a contradiction.\kraj

For the rest of this section we need to work with a nonstandard extension which is not only ${\goth c}^+$-saturated, but in which all the sets $\{x\in\mathbb{N}:x\leq z\}$ for $z\in\zve{\mathbb{N}}\setminus\mathbb{N}$ also have the same cardinality. We denote
$$\infty=|\{x\in\mathbb{N}:x\leq z\}|.$$
By \cite{DN3}, Corollary 2.3, this condition is true in extensions satisfying a forcing-like axiom $\Delta_1$, introduced there. Models satisfying both ${\goth c}^+$-saturation and $\Delta_1$ can be constructed in ZFC; let us denote this joint assumption with $\Delta_1$+SAT.

\begin{de}
For $u,v\in{\cal E}_{\cal P}$ and $x\in\mathbb{N}$, we denote $D^{[u,v]}_x:=\{(p,k):u\preceq_{\cal P}[(p,k)]_{\approx_{\cal P}}\preceq_{\cal P}v\land p^{k}\zvepar x\}$.
\end{de}

Under $\Delta_1$+SAT, for any basic class ${\cal P}^u$, every generator of any ultrafilter $\cal F$ has the same amount of ingredients from $\mu({\cal P}^u)$, and the only possible infinite "quantity" of such ingredients is $\infty$. In fact, by \cite{So8}, Theorem 3.2, each set $D^{[u,v]}_x$ is either finite or has cardinality $\infty$.

This allows the following generalization of patterns; another change that happens when we move outside $L$ (beside adding higher powers of $\cal P$) is that ultrafilters can be divisible by the same basic class infinitely many times. 

\begin{de}\label{pattern}
Denote $\mathbb{N}_\infty=\omega\cup\{\infty\}$. Let ${\cal A}$ be the set of all functions $\alpha:{\cal B}\rightarrow\mathbb{N}_\infty$. Elements $\alpha\in{\cal A}$ are called patterns.

For any $x=\prod_{n\leq z}p_n^{h(n)}\in\zve{\mathbb{N}}$ as in Proposition \ref{fund}, define $\alpha_x\in{\cal A}$ as follows. For each basic ${\cal P}^u\in{\cal B}$, let $\alpha_x({\cal P}^u)=|D^{[u,u]}_x|$.
\end{de}

We introduce a preorder on the family of all patterns.

\begin{de}\label{defdominate}
Let $(L,\leq)$ be a linear order and let $a=\langle a_m:m\in L\rangle$, $b=\langle b_m:m\in L\rangle$ be two sequences in $\mathbb{N}_\infty$. We say that $a$ dominates $b$ if, for every $l\in L$:
$$\sum_{m\geq l}a_m\geq\sum_{m\geq l}b_m.$$

For $\alpha,\beta\in{\cal A}$ we define: $\alpha\preceq\beta$ if $\langle\beta({\cal P}^u):u\in{\cal E}_{\cal P}\rangle$ dominates $\langle\alpha({\cal P}^u):u\in{\cal E}_{\cal P}\rangle$ for every ${\cal P}\in\overline{\mathbb{P}}$. If $\alpha\preceq\beta$ and $\beta\preceq\alpha$, we write $\alpha\approx\beta$.
\end{de}

\begin{te}[\cite{So8}, Lemma 3.10]\label{pomveza}
Assume $\Delta_1$+SAT. If $x\zvemid y$, then $\alpha_x\preceq\alpha_y$.
\end{te}

\begin{te}[\cite{So8}, Theorem 3.11]\label{samealpha}
Assume $\Delta_1$+SAT. For any ${\cal F}\in\beta\mathbb{N}$ and any two $x,y\in\mu({\cal F})$ holds $\alpha_x=\alpha_y$.
\end{te}

Thus, under the assumption $\Delta_1$+SAT, just like in $L$ one can define $\alpha_{\cal F}:=\alpha_x$ for ${\cal F}\in\beta\mathbb{N}$ and any $x\models{\cal F}$. The following corollary ensues.

\begin{co}\label{vezaalphamid}
Assume $\Delta_1$+SAT.
 
(a) If ${\cal F}\widemid{\cal G}$, then $\alpha_{\cal F}\preceq\alpha_{\cal G}$.

(b) If ${\cal F}=_\sim{\cal G}$, then $\alpha_{\cal F}\approx\alpha_{\cal G}$.
\end{co}

To generalize Definition \ref{defFalpha} of the family of sets that must belong to all ultrafilters of a given pattern, we need to consider a
larger collection of upward-closed sets.

\begin{de}\label{defF}
Let $\alpha\in{\cal A}$, ${\cal P}\in\overline{\mathbb{P}}$ and $A\in{\cal P}\upharpoonright\mathbb{P}$.

For any $h:A\rightarrow\mathbb{N}$ let
$$A^h=\{m\in\mathbb{N}:(\exists p\in A)p^{h(p)}\mid m\}.$$

For $u\in{\cal E}_{\cal P}$, an $(A,{\cal P}^{u})$-set is any set of the form $A^h$ such that for some/every $(p,x)\in u$ holds $p^x\in\zve A^h$.

An $(A,{\cal P}^{w})$-set for some $w\succeq_{\cal P}u$ which is not an $(A,{\cal P}^{v})$-set for any $v\prec_{\cal P}u$ will be called an $(A,{\cal P}^{\succeq u})$-set.

An $(\alpha,A,{\cal P})$-set is any finite product of $(A,{\cal P}^u)$-sets for various $u\in{\cal E}_{\cal P}$, such that for any fixed $u$, if $\sum_{w\succeq_{\cal P} u}\alpha({\cal P}^w)=n\in\mathbb{N}$, then there are at most $n$ $(A,{\cal P}^{\succeq u})$-sets in the product.

An $\alpha$-set is any finite product $C_1C_2\dots C_k$ of $(\alpha,A_i,{\cal P}_i)$-sets $C_i$, with $A_i\in{\cal P}_i$, ${\cal P}_i\neq{\cal P}_j$ and $A_i\cap A_j=\emptyset$ for $i\neq j$.

Finally, $F_\alpha$ is the intersection of ${\cal U}$ with the filter generated by the family of all $\alpha$-sets.
\end{de}

Now, if $\sigma(\alpha)\in\mathbb{N}$ (see Definition \ref{nivoi}), for any set $B$ belonging to $F_\alpha$ in the sense of Definition \ref{defFalpha}, it is easy to see that $B\uparrow$ is an $F_\alpha$-set in the sense of Definition \ref{defF}.

\begin{te}[\cite{So8}, Theorem 4.3]\label{efovi}
For every ${\cal F}\in\beta\mathbb{N}$, $F_{\alpha_{\cal F}}\subseteq {\cal F}\cap{\cal U}$.
\end{te}

In fact, if we narrow $\cal U$ to the family ${\cal U}'$ consisting only of sets of the form $A^h$, then $F_{\alpha_{\cal F}}={\cal F}\cap{\cal U}'$; see \cite{DLMS} for details. However, the inclusion in the previous theorem can be strict; in other words the pattern $\alpha_{\cal F}$ does not necessarily determine the $=_\sim$-equivalence class of $\cal F$. For example, if ${\cal P}\in\overline{\mathbb{P}}$ is not Ramsey, by Theorem \ref{ramsey} there are $=_\sim$-nonequivalent ultrafilters with the same pattern $\{({\cal P},2)\}$.

The pattern of an ultrafilter can not be just any function from ${\cal B}$ to $\mathbb{N}_\infty$; it must satisfy an additional condition that we call $\cal U$-closedness. Namely, on $\cal B$ we can define a topology with $\overline{A}$ (for $A\in{\cal U}$) as base sets. (Note that, if ${\cal F},{\cal G}\in{\cal P}^u$, then $A\in{\cal F}$ holds if and only if $A\in{\cal G}$.

\begin{de}
A pattern $\alpha$ is $\cal U$-closed if, whenever $\sum_{w\succeq_{\cal Q}u}\alpha({\cal Q}^w)$ is finite, then there is a neighborhood $\overline{A}$ of ${\cal Q}^u$ in which there are no basic classes ${\cal P}^v$ such that $\alpha({\cal P}^v)>0$ other than ${\cal Q}^w$ for $w\succeq_{\cal Q}u$.

The family of all $\cal U$-closed patterns is denoted ${\cal A}_{cl}$.
\end{de}

\begin{te}[\cite{So8}, Theorem 3.9]\label{zatvorenost}
For every ${\cal F}\in\beta\mathbb{N}$, $\alpha_{\cal F}\in{\cal A}_{cl}$.
\end{te}

\begin{ex}
Recall that with $MAX$ we denoted the $\widemid$-greatest class. By \cite{So4}, Lemma 4.6, ${\cal F}\in MAX$ if and only if $m\widemid{\cal F}$ for all $m\in{\mathbb{N}}$. Hence, ${\cal F}\in MAX$ if and only if its pattern $\alpha_{\cal F}$ contains $\{(p^\omega,1):p\in\mathbb{P}\}$. By $\cal U$-closedness it follows that actually $\alpha_{\cal F}({\cal P}^{max})=\infty$ for every ${\cal P}\in\overline{P}$.

This enables us to conclude that Corollary \ref{vezaalphamid} can not be strengthened: ${\cal F}=_\sim{\cal G}$ does not imply $\alpha_{\cal F}=\alpha_{\cal G}$. Namely, take any $x\in\mu(MAX)$. Define $a=\prod_{p\in\zve{\mathbb{P}},p\zve{\mid}x}p$; this is well-defined by Proposition \ref{fund}. Now let $y=ax$ and, for any prime $p$, $z=py$. Clearly, $y,z\in\mu(MAX)$ as well, since they are divisible by $x$. However, if ${\cal P}=tp(p/\zve{\mathbb{N}})$, then $\alpha_y({\cal P})=0$ and $\alpha_z({\cal P})=1$.
\end{ex}

\begin{te}[\cite{So8}, Theorem 4.4]\label{vezaalpha}
For patterns $\alpha,\beta\in{\cal A}_{cl}$, the following conditions are equivalent:

(i) $\alpha\preceq\beta$;

(ii) $F_{\alpha}\subseteq F_{\beta}$.
\end{te}

\begin{te}[\cite{So8}, Theorem 4.7]\label{overunder}
Let $\beta\in{\cal A}_{cl}$ and ${\cal F}\in\beta\mathbb{N}$.

(a) If $\alpha_{\cal F}\preceq\beta$, then there is ${\cal G}\in\beta\mathbb{N}$ such that $\alpha_{\cal G}\approx\beta$ and ${\cal F}\widemid{\cal G}$.

(b) If $\beta\preceq\alpha_{\cal F}$, then there is ${\cal G}\in\beta\mathbb{N}$ such that $\alpha_{\cal G}\approx\beta$ and ${\cal G}\widemid{\cal F}$.
\end{te}

\section{Congruence modulo an ultrafilter}\label{sec6}

After considering divisibility, the next natural relation to try to extend to ultrafilters is congruence. For this, instead of $\beta\mathbb{N}$, we work with $\beta\mathbb{Z}=\beta\mathbb{N}\cup\beta(-\mathbb{N})\cup\{0\}$, where $-A=\{-a:a\in A\}$, $-{\cal F}=\{-A:A\in{\cal F}\}$ and $\beta(-\mathbb{N})=\{-{\cal F}:{\cal F}\in\beta\mathbb{N}\}$. On $\beta\mathbb{Z}$ we can, beside addition and multiplication, also use subtraction, defined either as ${\cal F}-{\cal G}={\cal F}+(-{\cal G})$ or as the extension of minus on $\mathbb{Z}$ using (\ref{extoper}); this can be done even though $-$ is not associative.

In fact, congruence of ultrafilters modulo $m\in\mathbb{N}$ was mentioned briefly in \cite{HS}. There, $\equiv_m$ is the kernel relation of the homomorphism $\widetilde{h_m}$ obtained as the extension of $h_m:\mathbb{Z}\rightarrow\{0,1,\dots,m-1\}$ defined by: $h_m(n)$ is the remainder of $n$ modulo $m$. Another way to do this is to use the same idea as in Definition \ref{defdiv}: ${\cal F}\equiv_m{\cal G}$ if and only if $(\forall A\in{\cal F})\{n\in\mathbb{Z}:(\exists a\in A)a\equiv_mn\}\in{\cal G}.$ Fortunatelly, these two definitions are equivalent. Several results on the relationship of $\equiv_m$ with $\widemid$ can be found in Section 2 of \cite{So6}.

For congurence modulo a nonprincipal ultrafilter, however, the task is not so easy. In \cite{So6} we proposed two extensions, both in accordance with the above definition of $\equiv_m$ for $m\in\mathbb{N}$. The first one again uses the idea from Definition \ref{defdiv}, this time not fixing $m$ but thinking of $x\equiv_my$ as a relation between $m$ and the pair $(x,y)$.

\begin{de}\label{defcong}
For ${\cal M}\in\beta\mathbb{N}$ and ${\cal F},{\cal G}\in\beta\mathbb{Z}$, ${\cal F}\equiv_{\cal M}{\cal G}$ if and only if for every $A\in{\cal M}$ the set $\{(x,y)\in\mathbb{Z}\times\mathbb{Z}:(\exists m\in A)x\equiv_my\}$ belongs to the ultrafilter ${\cal F}\otimes{\cal G}$.
\end{de}

\begin{te}[\cite{So6}, Lemma 4.2 and Theorem 4.5]\label{equivdelj}
For ${\cal M}\in\beta\mathbb{N}$ and ${\cal F},{\cal G}\in\beta\mathbb{Z}$, the following conditions are equivalent:

(i) ${\cal F}\equiv_{\cal M}{\cal G}$;

(ii) ${\cal M}\widemid{\cal F}-{\cal G}$;

(iii) in some ${\goth c}^+$-saturated nonstandard extension holds
\begin{equation}\label{eqpair}
(\forall m\in\mu({\cal M}))(\exists x\in\mu({\cal F}))(\exists y\in\mu({\cal G}))((x,y)\mbox{ is a tensor pair }\land m\zvemid x-y)
\end{equation}

(iv) in every ${\goth c}^+$-saturated nonstandard extension holds (\ref{eqpair}).
\end{te}

Unfortunately, $\equiv_{\cal M}$ is not an equivalence relation for all ${\cal M}$. It is reflexive, but it is shown in \cite{DLM} (Corollary 2.6 and Theorem 2.3) that it need not be symmetric and transitive.

The second type of extension, called strong congruence, originally used iterated hyperextensions. In \cite{DLM} an equivalent condition was proposed that does not require such advanced machinery.

\begin{de}\label{defstrong}
For ${\cal M}\in\beta\mathbb{N}$ and ${\cal F},{\cal G}\in\beta\mathbb{Z}$, ${\cal F}\equiv^s_{\cal M}{\cal G}$ if and only if 
$$(\exists m\in\mu({\cal M}))(\exists x\in\mu({\cal F}))(\exists y\in\mu({\cal G}))((m,x,y)\mbox{ is a tensor triple }\land m\zvemid x-y).$$
\end{de}

What makes this definition convenient is the following observation.

\begin{pp}
${\cal F}\equiv^s_{\cal M}{\cal G}$ implies
$$(\forall m\in\mu({\cal M}))(\forall x\in\mu({\cal F}))(\forall y\in\mu({\cal G}))((m,x,y)\mbox{ is a tensor triple }\Rightarrow m\zvemid x-y).$$
\end{pp}

\dokaz Let $m\in\mu({\cal M})$, $x\in\mu({\cal F})$ and $y\in\mu({\cal G})$ be such that $m\zvemid x-y$. This means that, if we define $A=\{(n,a,b)\in\mathbb{N}\times\mathbb{N}\times\mathbb{N}:n\mid a-b\}$, then $(m,x,y)\in\zve A$. If $(m,x,y)$ is a tensor triple, then $(m,x,y)\models{\cal M}\otimes{\cal F}\otimes{\cal G}$, so it follows that $A\in{\cal M}\otimes{\cal F}\otimes{\cal G}$. Now if $m'\in\mu({\cal M})$, $x'\in\mu({\cal F})$ and $y'\in\mu({\cal G})$ are arbitrary such that $(m',x',y')$ is a tensor triple, we get $(m',x',y')\in\zve A$ as well.\kraj

This relation is an equivalence relation, in fact a congruence with respect to both addition and multiplication.

\begin{te}[\cite{So6}, Theorem 5.7]\label{compat}
Let ${\cal M}\in\beta\mathbb{N}$. $\equiv_{\cal M}^s$ is compatible with operations $+$ and $\cdot$ on $\beta\mathbb{Z}$:

(a) ${\cal F}_1\equiv_{\cal M}^s{\cal F}_2$ and ${\cal G}_1\equiv_{\cal M}^s{\cal G}_2$ imply ${\cal F}_1+{\cal G}_1\equiv_{\cal M}^s{\cal F}_2+{\cal G}_2$;

(b) ${\cal F}_1\equiv_{\cal M}^s{\cal F}_2$ and ${\cal G}_1\equiv_{\cal M}^s{\cal G}_2$ imply ${\cal F}_1\cdot{\cal G}_1\equiv_{\cal M}^s{\cal F}_2\cdot{\cal G}_2$.
\end{te}

On the other hand, this relation has another drawback: ${\cal M}\equiv^s_{\cal M}0$ is not true for all ultrafilters $\cal M$. 

\begin{de}
Let ${\cal F},{\cal G}\in\beta\mathbb{N}$. The strong divisibility relation is defined by: ${\cal F}\mid^s{\cal G}$ if ${\cal G}\equiv_{\cal F}0$.

Le also $D({\cal G})=\{n\in\mathbb{N}:n\mid{\cal G}\}$.
\end{de}

\begin{lm}[\cite{DLM}, Remark 1.4 and Theorem 3.10]
The following conditions are equivalent:

(i) ${\cal F}\mid^s{\cal G}$;

(ii) $(\exists x\in\mu({\cal F})(\exists y\in\mu({\cal G})((x,y)\mbox{ is a tensor pair }\land x\zvemid y)$;

(iii) $(\forall x\in\mu({\cal F})(\forall y\in\mu({\cal G})((x,y)\mbox{ is a tensor pair }\Rightarrow x\zvemid y)$;

(iv) $D({\cal G})\in{\cal F}$.
\end{lm}

Thus, the strong divisibility relation is too strong: not only it is not reflexive, but it only depends on divisors from $\mathbb{N}$, so $\mathbb{N}$-free ultrafilters do not have any strong divisors at all.

To make any of the two introduced congruence relations useful, it is important to establish for which ultrafilters $\cal M$ they behave nicely. It turned out that in both cases the answer to this problem is the same.

\begin{de}
${\cal F}$ is self-divisible if $D({\cal F})\in{\cal F}$.
\end{de}

Thus, $\cal F$ is self-divisible if ${\cal F}\equiv^s_{\cal F}0$. 

\begin{te}[\cite{DLM}, Theorem 3.10]
For any ${\cal F}\in\beta\mathbb{N}$ the following conditions are equivalent:

(1) $\cal F$ is self-divisible;

(2) the relations $\equiv_{\cal F}$ and $\equiv_{\cal F}^s$ coincide;

(3) $\equiv_{\cal F}$ is an equivalence relation;

(4) for every $a,b\models{\cal F}$ there is $c\models{\cal F}$ such that $c\zvemid gcd(a,b)$.
\end{te}

Many more equivalent conditions are listed in \cite{DLM}, Theorem 6.7. Another one using patterns (as in Definition \ref{pattern}) is due to appear in \cite{DLMS}.

\section{Other divisibility-like relations}\label{sec7}

$\widemid$ was only one of several extensions of the divisibility proposed in \cite{So1}. Although others do not have so many properties reflecting the properties of $\mid$ on $\mathbb{N}$, they are interesting in other ways.

\begin{de}
Let ${\cal F},{\cal G},{\cal H}\in\beta\mathbb{N}$.

(a) $\cal G$ is left-divisible by ${\cal F}$ (${\cal F}\mid_L{\cal G}$) if there is ${\cal H}\in\beta\mathbb{N}$ such that ${\cal G}={\cal H}\cdot{\cal F}$.

(b) ${\cal G}$ is right-divisible by ${\cal F}$ (${\cal F}\mid_R{\cal G}$) if there is ${\cal H}\in\beta\mathbb{N}$ such that ${\cal G}={\cal F}\cdot{\cal H}$.

(c) ${\cal G}$ is mid-divisible by ${\cal F}$ (${\cal F}\mid_M{\cal G}$) if there are ${\cal H},{\cal K}\in\beta\mathbb{N}$ such that ${\cal G}={\cal H}\cdot{\cal F}\cdot{\cal K}$.
\end{de}

In semigroup theory $\mid_L$, $\mid_R$ and $\mid_M$ are usually called Green preorders; for example ${\cal F}\mid_L{\cal G}$ means that the principal left ideal ${\cal G}\cdot\beta\mathbb{N}$ is included in ${\cal F}\cdot\beta\mathbb{N}$. $\mid_R$ is the analogue for principal right ideals, and $\mid_M$, the transitive closure of $\mid_L\cup\mid_R$, the analogue for two-sided ideals.

It was shown in \cite{So1} that $\widemid$ is the maximal continuous extension of divisibility, meaning that the "pre-image" $\{{\cal G}\in\beta\mathbb{N}:{\cal G}\widemid{\cal F}\}$ of any singleton $\{{\cal F}\}$ is closed. $\mid_R$ was also proved to be continuous.

There is another preorder on $\beta\mathbb{N}$ which is, strictly speaking, not a divisibility relation, but is nicely intertwined between these relations.

\begin{de}
For $A,B\in P({\mathbb{N}})$, $A$ is finitely embeddable in $B$ ($A\leq_{fe}B$) if for every finite $F\subseteq A$ there is $k\in {\mathbb{N}}$ such that $kF\subseteq B$.

For ${\cal F},{\cal G}\in\beta {\mathbb{N}}$, ${\cal F}$ is finitely embeddable in ${\cal G}$ (${\cal F}\leq_{fe}{\cal G}$) if for every $B\in{\cal G}$ there is $A\in{\cal F}$ such that $A\leq_{fe}B$.
\end{de}

A similar relation, based on addition instead of multiplication, was considered earlier in \cite{BD}, \cite{DN1} and \cite{L2}. A general study of finite embeddability relations was carried out by Luperi Baglini in \cite{L3}. It turned out that this relation lays between some of the mentioned divisibility relations.

\begin{te}[\cite{So7}, Theorem 28]\label{feconnection}
$$\begin{array}{c}
\mid_L\\
\\
\mid_R
\end{array}\hspace{-1mm}
\begin{array}{c}
\rotatebox[origin=c]{-45}{$\subset$}\\
\rotatebox[origin=c]{45}{$\subset$}\\
\end{array}\hspace{-1mm}
\;\mid_M\;\subset\;\leq_{fe}\;\subset\hspace{1mm}\widetilde{\mid}\hspace{1mm}$$
\end{te}

Also, by \cite{So7}, Lemma 2.3, $\leq_{fe}$ coincides with each of the divisibility relations when divisibility by $m\in\mathbb{N}$ is concerned. The connection from Theorem \ref{feconnection} makes it possible to better understand the $\leq_{fe}$-hierarchy by considering the (weaker) $\widemid$-hierarchy, and vice versa. For more about this connection we address the reader to \cite{So7}.

\section{Remarks and open problems}

Here we collect several problems which remained unsolved so far, and make some final remarks.

The organization of $L$ into levels makes it much easier to work with. The order $(\zve{\mathbb{N}},\zvemid)$ is also organized into levels, which can be seen by considering the extension $\zve lev$ of the function $lev:\zve{\mathbb{N}}\rightarrow\zve{\mathbb{N}}$ defined by $lev(p_1^{a_1}p_2^{a_2}\dots p_k^{a_k})=a_1+a_2+\dots+a_k$, see \cite{So4}, Lemma 2.6. If one tried to do the same for ultrafilters, a function $\widetilde{lev}:\beta\mathbb{N}\rightarrow\beta\mathbb{N}$ would be obtained. However, as mentioned in Section \ref{sec5}, the orders $({\cal E}_{\cal P},\preceq_{\cal P})$ need not be isomorphic for different $\cal P$, and they would have to be if their members belonged to corresponding levels of the $\widemid$-hierarchy. Hence it is not possible to define a system of levels in the whole $\beta\mathbb{N}$, not even of non-well-ordered type.

The question: which ultrafilters have immediate predecessors, and which can be written as limits of $\widemid$-increasing chains, will be addressed in \cite{DLMS}. A general result on the exact number of immediate predecessors or successors still seems to be out of reach.

Another unresolved problem is Question 5.2 from \cite{So8}: is it possible to improve Corollary \ref{overunder} to get an ultrafilter $\cal G$ such that $\alpha_{\cal G}=\beta$ (and not only $\alpha_{\cal G}\approx\beta$)?\\

This research was supported by the Science Fund of the Republic of Serbia, Grant No. 7750027 (Set-theoretic, model-theoretic and Ramsey-theoretic phenomena in mathematical structures: similarity and diversity -- SMART), and by the Ministry of Science, Technological Development and Innovation of the Republic of Serbia (grant no.\ 451-03-47/2023-01/200125).

\end{document}